\RequirePackage{ifpdf}
\ifpdf 
\documentclass[pdftex]{sigma}
\else
\documentclass{sigma}
\fi

\numberwithin{equation}{section}

\begin{document}

\allowdisplaybreaks

\renewcommand{\thefootnote}{$\star$}

\renewcommand{\PaperNumber}{092}

\FirstPageHeading

\ShortArticleName{An Introduction to the $q$-Laguerre--Hahn Orthogonal $q$-Polynomials}

\ArticleName{An Introduction to the $\boldsymbol{q}$-Laguerre--Hahn\\ Orthogonal
$\boldsymbol{q}$-Polynomials\footnote{This
paper is a contribution to the Proceedings of the Conference ``Symmetries and Integrability of Dif\/ference Equations (SIDE-9)'' (June 14--18, 2010, Varna, Bulgaria). The full collection is available at \href{http://www.emis.de/journals/SIGMA/SIDE-9.html}{http://www.emis.de/journals/SIGMA/SIDE-9.html}}}

\Author{Abdallah GHRESSI, Lotfi KH\'{E}RIJI and Mohamed Ihsen TOUNSI}

\AuthorNameForHeading{A.~Ghressi, L.~Kh\'{e}riji and M.I.~Tounsi}

\Address{Institut Sup\'{e}rieur des Sciences Appliqu\'{e}es et de
Technologies de Gab\`{e}s,\\ Rue Omar Ibn El-Khattab 6072 Gab\`{e}s,
Tunisia}
\Email{\href{mailto:Abdallah.Ghrissi@fsg.rnu.tn}{Abdallah.Ghrissi@fsg.rnu.tn}, \href{mailto:kheriji@yahoo.fr}{kheriji@yahoo.fr}, \href{mailto:MohamedIhssen.Tounsi@issatgb.rnu.tn}{MohamedIhssen.Tounsi@issatgb.rnu.tn}}

\ArticleDates{Received February 14, 2011, in f\/inal form September 26, 2011;  Published online October 04, 2011}

\Abstract{Orthogonal $q$-polynomials associated with
$q$-Laguerre--Hahn form will be studied as a generalization of the
$q$-semiclassical forms via a suitable $q$-dif\/ference equation. The
concept of class and a criterion to determinate it will be given.
The $q$-Riccati equation satisf\/ied by the corresponding formal
Stieltjes series is obtained. Also, the structure relation is
established. Some illustrative examples are highlighted.}

\Keywords{orthogonal $q$-polynomials; $q$-Laguerre--Hahn form;
$q$-dif\/fe\-rence operator; $q$-dif\/ference equation; $q$-Riccati
equation}

\Classification{42C05; 33C45}

\section{Introduction and preliminary results}

The concept of the usual Laguerre--Hahn polynomials were extensively
studied by several authors \cite{1,2,4,6,8,9,10,15,18}. They constitute a
very remarkable family of orthogonal polynomials taking
consideration of most of the monic orthogonal polynomials
sequences (MOPS) found in lite\-rature. In
particular, semiclassical orthogonal polynomials are Laguerre--Hahn
MOPS~\cite{15,20}. The Laguerre--Hahn set of form (linear functional) is
invariant under the standard perturbations of forms \cite{2,9,18,20}. It
is well known that a usual Laguerre--Hahn polynomial satisf\/ies a
fourth order dif\/ferential equation with polynomials coef\/f\/icients but
the converse remains not proved until now~\cite{20}. Discrete
Laguerre--Hahn polynomials were studied in~\cite{13}. These families are
already extensions of discrete semiclassical polynomials~\cite{19}. In
literature, analysis and characterization of the $q$-Laguerre--Hahn
orthogonal $q$-polynomials have not been yet presented in a unif\/ied
way. However, several authors have studied the fourth order
$q$-dif\/ference equation related to some examples of
$q$-Laguerre--Hahn orthogonal $q$-polynomials such as the
co-recursive and the $r$th associated of $q$-classical
polynomials~\cite{11,12}. More generally, the fourth order dif\/ference
equation of Laguerre--Hahn orthogonal on special non-uniform lattices
polynomials was established in~\cite{4}. For other relevant works in the
domain of orthogonal $q$-polynomials and $q$-dif\/ference
equation theory see~\cite{3,21} and~\cite{5}.

So the aim of this contribution is to establish a basic theory of
$q$-Laguerre--Hahn orthogonal $q$-polynomials. We give some
characterization theorems for this case such as the structure
relation and the $q$-Riccati equation. We extend the concept of the
class of the usual Laguerre--Hahn forms to the $q$-Laguerre--Hahn
case. Moreover, we show that some standard transformation and
perturbation carried out on the $q$-Laguerre--Hahn forms lead to new
$q$-Laguerre--Hahn forms; the class of the resulting forms is
analyzed and some examples are treated.

We denote by $\mathcal{P}$ the vector space of the polynomials with
coef\/f\/icients in $\mathbb{C}$ and by $\mathcal{P}'$ its dual space
whose elements are forms. The action of $u\in \mathcal{P}^{\prime }$
on $f\in \mathcal{P}$ is denoted as $\langle u,f\rangle $. In
particular, we denote by $(u)_{n}:=\langle u,x^{n}\rangle
\hspace{0.1cm},\hspace{0.1cm}n\geq 0$ the moments of $u$. A linear
operator $T : \mathcal{P}\longrightarrow \mathcal{P}$ has a
transpose $^{t}T : \mathcal{P}'\longrightarrow \mathcal{P}'$ def\/ined
by
\[
\langle ^{t}T u,f\rangle=\langle u,Tf\rangle , \qquad u\in
\mathcal{P}^{\prime }, \qquad f\in \mathcal{P}.
\]
For instance, for any form $u$, any polynomial $g$ and any
$(a,c)\in(\mathbb{C}\setminus \{0\})\times\mathbb{C}$, we let
$H_{q}u$,   $gu$,  $h_{a}u$,  $D u$,  $(x-c)^{-1}u$ and~$\delta_{c}$, be the forms def\/ined as usually \cite{20} and~\cite{16} for the
results related to the operator $H_{q}$
\begin{gather*}
\langle H_{q}u,f \rangle :=-\langle u,H_{q}f\rangle
 ,\qquad  \langle gu,f\rangle:=\langle
u,gf\rangle ,\qquad \langle
h_{a}u,f\rangle:=\langle u,h_{a}f\rangle ,
\\
\langle D u,f\rangle:=-\langle
u,f'\rangle ,\qquad \langle (x-c)^{-1}u,f\rangle
:=\langle u,\theta _{c}f\rangle ,\qquad \langle
\delta_{c},f\rangle:=f(c) ,
\end{gather*}
where for all $f\in \mathcal{P}$ and $q\in
\widetilde{\mathbb{C}}:=\big\{z\in \mathbb{C}, \ z\neq
0, \ z^{n}\neq 1, \ n \geq 1\big\}$ \cite{16}
\[
(H_{q}f)(x)=\frac{f(qx)-f(x)}{(q-1)x} ,\qquad
(h_{a}f)(x)=f(ax) ,\qquad (\theta_{c}f)(x)
=\frac{f(x) -f(c)}{x-c}.
\]
In particular, this yields to
\[
(H_{q}u)_{n}=-[n]_{q}(u)_{n-1} ,\qquad n\geq 0 ,
\]
 where $(u)_{-1}=0$ and
$[n]_{q}:= \frac{q^{n}-1}{q-1}$, $n\geq
0$~\cite{15}.  It is obvious that when $q \rightarrow 1$, we
meet again the derivative $D$.

For $f\in \mathcal{P}$ and $u\in \mathcal{P}^{\prime }$, the product
$uf$  is the polynomial~\cite{20}
\[
(uf)(x):=\langle u, \frac{xf(x)-\zeta f(\zeta) }{x-\zeta
}\rangle=\sum_{i=0}^{n}\left(\sum_{j=i}^{n}(u)_{j-i}\, f_{j}\right)
x^{i},
\]
where $f(x)= \sum\limits_{i=0}^{n} f_{i} x^{i}$. This allows us
to def\/ine the Cauchy's product of two forms:
\[
\langle uv,f\rangle:=\langle
u,vf\rangle ,\qquad f\in \mathcal{P}.
\]
The product def\/ined as before is commutative~\cite{20}. Particularly, the
inverse $u^{-1}$ of $u$ if there exists is def\/ined by $u   u^{-1}=\delta_{0}$.

The Stieltjes formal series of $u\in \mathcal{P}^{\prime }$ is
def\/ined by
\[
S(u)(z) :=- \sum_{n\geq 0}
 \frac{(u)_{n}}{z^{n+1}}.
\]

A form $u$ is said to be regular whenever there is a sequence
 of monic polynomials $\{P_{n}\}_{n \geq 0}$, $\deg P_{n}=n$,
$n \geq0$  such that $\langle u,P_{n} P_{m}\rangle=r_{n}
\delta_{n,m}$ with $r_{n}\neq 0$ for any $n, m \geq 0$. In this
case, $\{P_{n}\}_{n \geq 0}$ is called a monic orthogonal
polynomials sequence MOPS and it is characterized by the following
three-term recurrence relation (Favard's theorem)
  \begin{gather}
   P_{0}(x)=1 ,\qquad  P_{1}(x)=x-\beta_{0} , \nonumber\\
   P_{n+2}(x)=(x- \beta_{n+1})P_{n+1}(x)-\gamma_{n+1} P_{n}(x),\qquad  n \geq 0 ,
  \label{1.1}
\end{gather}
where $\beta_{n}= \frac{\langle u,x P_{n}^{2}
\rangle}{r_{n}} \in \mathbb{C}$,
$\gamma_{n+1}= \frac{r_{n+1}}{r_{n}} \in
\mathbb{C}\setminus\{0\}$, $ n \geq 0$.

The shifted MOPS
$\{\widehat{P}_{n}:=a^{-n} (h_{a}P_{n} )\}_{n\geq0}$ is then
orthogonal with respect to $\widehat{u}=h_{a^{-1}}u$ and satisf\/ies
\eqref{1.1} with~\cite{20}
\[
\widehat{\beta}_{n}=\frac{\beta_{n}}{a}   , \qquad
\widehat{\gamma}_{n+1}=\frac{\gamma_{n+1}}{a^{2}}  , \qquad
n\geq0.
\]
Moreover, the form $u$ is said to be normalized if $(u)_{0}=1$. In
this paper, we suppose that any form will be normalized.

The form $u$ is said to be positive def\/inite if and only if
$\beta_{n} \in \mathbb{R}$ and $\gamma_{n+1}>0$ for all $n \geq 0$.
When $u$ is regular, $\{P_{n}\}_{n\geq 0}$ is a symmetrical MOPS if
and only if $\beta_{n}=0$, $n\geq0$ or equivalently
$(u)_{2n+1}=0$, $n\geq0$.

Given a regular form $u$ and the corresponding MOPS
$\{P_n\}_{n\geq0}$, we def\/ine the associated sequence of the f\/irst
kind $\big\{P_{n}^{(1)}\big\}_{n\geq0}$ of $\{P_n\}_{n\geq0}$ by \cite[equations~(2.8) and~(2.9)]{20}
\[
{P_n}^{(1)}(x)= \bigl\langle u,
\frac{P_{n+1}(x)-P_{n+1}(\xi)}{x-\xi}\bigr\rangle=(u\theta_0
P_{n+1})(x), \qquad n\geq0.
\]

The following well known results (see~\cite{16,17,20}) will be
needed in the sequel.
\begin{lemma}\label{lemma1} Let $u \in \mathcal{P}'$. $u$ is regular if and only if
 $\Delta_{n}(u)\neq 0$,
$n\geq 0$ where
\[
\Delta_{n}(u):=\det\big((u)_{\mu+\nu}\big)_{\mu,\nu=0}^{n} ,
\qquad n\geq 0
\]
are the Hankel determinants.
\end{lemma}

\begin{lemma}\label{lemma2} For  $f,g\in \mathcal{P}$, $u, v \in \mathcal{P}^{\prime }$,
$(a,b,c)\in \mathbb{C}\setminus\{0\}\times\mathbb{C}^{2}$, and
$n\geq1$, we have
\begin{gather}
(x-c)\big((x-c)^{-1}u\big) =u   , \qquad (x-c)^{-1}((x-c) u) =u-(u)
_{0}\delta _{c}   , \label{1.2}\\
(u\theta_{0}f)(x)=a_{n} x^{n-1}(u)_{0}+\textrm{lower order terms}
,  \qquad f(x)=\sum_{k=0}^{n}a_{k}x^{k},\label{1.3}
\\
u\theta_0 (fg)=g(u\theta_0 f)+(fu)\theta_0 g   , \label{1.4}
\\
u\theta_{0}(f P_{k+1})=f P_{k}^{(1)}, \qquad k+1\geq \deg f,
\label{1.5}
\\
\theta_{b}-\theta_{c}=(b-c) \theta_{b}\circ\theta_{c}  , \qquad
\theta_{b}\circ\theta_{c}=\theta_{c}\circ\theta_{b} ,\label{1.6}
\\
h_{a}(gu) =( h_{a^{-1}}g) (h_{a}u)   , \qquad h_{a}(uv) =( h_{a}u)
(h_{a}v)   , \qquad h_{a}\big(x^{-1} u\big)= a   x^{-1} h_{a}u, \label{1.7}
\\
h_{q^{-1}}\circ H_{q}=H_{q^{-1}}   , \qquad H_{q}\circ
h_{q^{-1}}=q^{-1}H_{q^{-1}}  ,\qquad \textit{in}\quad  \mathcal{P
},\label{1.8}
\\
h_{q^{-1}}\circ H_{q}=q^{-1}H_{q^{-1}}   ,\qquad H_{q}\circ
h_{q^{-1}}=H_{q^{-1}} ,
\qquad \textit{in}\quad \mathcal{P}^{\prime
},\label{1.9}
\\
H_{q}(fg)(x) =(h_{q}f)(x)(H_{q}g)(x) +g(x)( H_{q}f)(x),\label{1.10}
\\
H_{q}(gu) =( h_{q^{-1}}g) H_{q}u+q^{-1}(H_{q^{-1}}g) u, \label{1.11}
\\
H_{q^{-1}}(u\theta _{0}f)(x) =q(H_{q}u)\theta _{0}(h_{q^{-1}}f)(x)
+( u\theta _{0}H_{q^{-1}}f)(x),\label{1.12}
\\
S(fu)(z) =f(z) S(u)(z) +(u\theta _{0}f)(z)   , \label{1.13}
\\
S(uv)(z) =-zS(u)(z)   S(v)(z) ,\label{1.14}
\\
S(x^{-n}u)(z) =z^{-n} S(u)(z)   , \qquad S(u^{-1})(z)=z^{-2}
 (S(u)(z) )^{-1} , \label{1.15}
\\
S(H_{q}u)(z) =q^{-1}(H_{q^{-1}}(S(u)))(z)   , \qquad
(h_{q^{-1}}S(u))(z)=q S(h_{q}u)(z).\label{1.16}
\end{gather}
\end{lemma}

\begin{definition}\label{definition1} A form $u$ is called
$q$-Laguerre--Hahn  when it is regular and satisf\/ies the
$q$-dif\/ference equation
\begin{gather}
 H_{q}(\Phi u) +\Psi u+B\big(x^{-1}u(h_{q}u)\big)=0,\label{1.17}
\end{gather}
where  $\Phi$, $\Psi$, $B$ are polynomials, with $\Phi$ monic. The
corresponding orthogonal sequence $\{ P_{n}\} _{n\geq 0}$ is called
$q$-Laguerre--Hahn MOPS.
\end{definition}

\begin{remark}\label{remark1} When
$B=0$ and the form $u$ is regular then $u$ is $q$-semiclassical
[17]. When $u$ is regular and not $q$-semiclassical then $u$ is
called a strict $q$-Laguerre--Hahn form.
\end{remark}

\begin{lemma}\label{lemma3} Let $u$ be a regular form. If $u$ is a strict $q$-Laguerre--Hahn form satisfying
\eqref{1.17} and there exist two polynomials $\Delta$ and $\Omega$ such
that
\begin{gather}
\Delta u+\Omega \big(x^{-1}u(h_{q}u)\big)=0 \label{1.18}
\end{gather}
then  $\Delta=\Omega=0$.
\end{lemma}

\begin{proof} The operation $\Delta \times \eqref{1.17}  - B \times \eqref{1.18}$ gives
\[
\Omega  H_{q}(\Phi u)+ (\Omega \Psi-\Delta B )u=0.
\]
According to \eqref{1.9} and \eqref{1.11}, the above equation becomes
\[
H_{q} ((h_{q}\Omega)\Phi u )+ (\Omega\Psi-
(H_{q}\Omega)\Phi-\Delta B )u=0.
\]
Then $\Delta=\Omega=0$ because the form $u$ is regular and not
$q$-semiclassical.
\end{proof}

\begin{lemma}\label{lemma4} Consider the
sequence $\{\widehat{P}_{n}\}_{n\geq0}$ obtained by shifting
$P_{n}$, i.e.\ $\widehat{P}_{n}(x)=a^{-n}P_{n}(ax)$, $n\geq0$, $a\neq0$.
When $u$ satisfies \eqref{1.17}, then $\widehat{u}=h_{a^{-1}}u$ fulfills
the $q$-difference equation
\[
H_{q}( \widehat{\Phi} \widehat{u}) +\widehat{\Psi}
\widehat{u}+\widehat{B}\big(x^{-1}\widehat{u}(h_{q}\widehat{u})\big)=0
 ,
\]
where
\hspace{0.1cm}$\widehat{\Phi}(x)=a^{-\deg\Phi}\Phi(ax)$,
$\widehat{\Psi}(x)=a^{1-\deg\Phi}\Psi(ax)$,
$\widehat{B}(x)=a^{-\deg\Phi}B(ax)$.
\end{lemma}

\begin{proof} With $u=h_{a}\widehat{u}$, we have
$\Psi u=\Psi
(h_{a}\widehat{u})=h_{a}\big((h_{a}\Psi)\widehat{u}\big)$
 from~\eqref{1.7}.
Further,
 \[
 H_{q}(\Phi u)=H_{q}\big(\Phi
(h_{a}\widehat{u})\big)=H_{q}\big(h_{a}\big((h_{a}\Phi)\widehat{u}\big)\big)
=a^{-1}h_{a}\big(H_{q}\big((h_{a}\Phi)\widehat{u}\big)\big)
\]
from \eqref{1.7} and \eqref{1.9}.

Moreover, by virtue of~\eqref{1.7} an other time we get
\[
B\big(x^{-1}u(h_{q}u)\big)=B\big(x^{-1}(h_{a}\widehat{u}) (h_{aq}
\widehat{u})\big)=B \big(x^{-1}h_{a}\big(\widehat{u}
h_{q}\widehat{u}\big)\big)
=a^{-1}h_{a}\big((h_{a}B)\big(x^{-1}\widehat{u}(h_{q}\widehat{u})\big)\big).
\]
Equation \eqref{1.17} becomes
\[
h_{a}\big(H_{q}\big(\Phi(ax)\widehat{u}\big)+
a\Psi(ax)\widehat{u}+B(ax)\big(x^{-1}\widehat{u}(h_{q}\widehat{u})\big)\big)=0.
\]
Hence the desired result.
\end{proof}

\section[Class of a $q$-Laguerre-Hahn form]{Class of a $\boldsymbol{q}$-Laguerre--Hahn form}

It is obvious that a $q$-Laguerre--Hahn form satisf\/ies an inf\/inite
number of $q$-dif\/ference equations type~\eqref{1.17}. Indeed, multiplying
\eqref{1.17} by a polynomial $\chi $ and taking into account~\eqref{1.7}, \eqref{1.11}
we obtain
\begin{gather}
H_{q}\big((h_{q}\chi)\Phi u\big)+\big\{\chi\Psi-\Phi
(H_{q}\chi)\big\}u+(\chi B)\big(x^{-1}u(h_{q}u)\big)=0.\label{2.1}
\end{gather}
Put $t=\deg \Phi$, $p=\deg
\Psi$, $r=\deg B$ with $d=\max(t,r)$ and
$s=\max(p-1,d-2)$. Thus, there exists $u\rightarrow \hbar(u)\subset
\mathbb{N}\cup\{-1\}$ from the set of $q$-Laguerre--Hahn forms into
the subsets of $\mathbb{N}\cup\{-1\}$.

\begin{definition}\label{definition2} The minimum element of $\hbar(u)$ will be called the
class of $u$. When $u$ is of class $s$, the sequence
$\{P_{n}\}_{n\geq0}$ orthogonal with respect to $u$ is said to be of
class $s$.
\end{definition}

\begin{proposition}\label{proposition1} The number $s$ is
an integer positive or zero. In other words, if $p=0$, then $d\geq
2$ or if    $0\leq d \leq 1$, then necessarily $p\geq 1$.
\end{proposition}

\begin{proof} Let us show that in case $s=-1$, the form
$u$ is not regular, which is a contradiction. Indeed, when $s=-1$,
we have
\[
\Phi(x)=c_{1}x+c_{0},\qquad \Psi(x)=a_{0}, \qquad
B(x)=b_{1}x+b_{0}
\]
with $c_{1}=1$ or $c_{1}=0$ and $c_{0}=1$, and where $a_{0}\neq0$.

The condition  $\langle H_{q}(\Phi u) +\Psi
u+B\big(x^{-1}u(h_{q}u)\big),x^{n}\rangle=0$,  $0\leq n
\leq 4$ gives successively
\begin{gather}
a_{0}+b_{1}=0,\nonumber\\
 (qb_{1}-c_{1} )(u)_{1}+b_{0}-c_{0}=0,
\label{2.2}
\\
\big(q^{2}b_{1}-(1+q)c_{1}\big)\big((u)_{2}-(u)_{1}^{2}\big)=0,
\label{2.3}
\\
\big(q^{3}b_{1}-\big(1+q+q^{2}\big)c_{1}\big)(u)_{3}+
\big\{\big(1+q^{2}\big)b_{0}+q(1+q)b_{1}(u)_{1}-\big(1+q+q^{2}\big)c_{0}\big\}(u)_{2}\nonumber\\
\qquad{}+qb_{0}(u)_{1}^{2}=0,\label{2.4}
\\
 \big(q^{4}b_{1}-(1+q)\big(1+q^{2}\big)c_{1}\big)(u)_{4}+
\big\{\big(1+q^{3}\big)b_{0}+q\big(1+q^{2}\big)b_{1}(u)_{1}-
(1+q)\big(1+q^{2}\big)c_{0}\big\}(u)_{3} \nonumber\\
\qquad{} +
q^{2}b_{1}(u)_{2}^{2}+q(1+q)b_{0}(u)_{1}(u)_{2}=0.\label{2.5}
\end{gather}
Suppose $q^{2}b_{1}-(1+q)c_{1}\neq 0.$ From \eqref{2.3}
\[
\Delta_{1}=\begin{vmatrix}
  1 & (u)_{1} \\
  (u)_{1} & (u)_{2}
\end{vmatrix}
=0.
\]
Contradiction.

Suppose $q^{2}b_{1}=(1+q)c_{1}=0$ implies $b_{1}=0=c_{1}$ implies
\eqref{2.2} $b_{0}=c_{0}=1$. Thus~\eqref{2.4} $(u)_{2}-(u)_{1}^{2}=0$, hence
$\Delta_{1}=0$. Contradiction.

Suppose $q^{2}b_{1}=(1+q)c_{1}\neq 0$ with $c_{1}=1$. From~\eqref{2.2} and
\eqref{2.4}, \eqref{2.5}, we have
\begin{gather}
(u)_{1}=q(c_{0}-b_{0}),\nonumber \\
(u)_{3} = q(c_{0}-2b_{0})(u)_{2}+q^{3}b_{0}(c_{0}-b_{0})^{2},\label{2.6}
\\
(u)_{4}=(u)_{2}^{2}+q^{2}b_{0}^{2}(u)_{2}-q^{4}b_{0}^{2}(c_{0}-b_{0})^{2}.\nonumber
\end{gather}
On the other hand, let us consider the Hankel determinant
\[
\Delta_{2}=\begin{vmatrix}
  1 & (u)_{1} & (u)_{2} \\
  (u)_{1} & (u)_{2} & (u)_{3} \\
  (u)_{2} & (u)_{3} & (u)_{4}
\end{vmatrix}.
\]
With \eqref{2.6}, we get $\Delta_{2}=0$. Contradiction.
\end{proof}

\begin{proposition}\label{proposition2}
Let $u$ be a strict $q$-Laguerre--Hahn form satisfying
\begin{gather}
  H_{q}(\Phi_1 u)+\Psi_1 u+B_1\big(x^{-1}u h_q u\big)=0
,\label{2.7}
\end{gather}
and
\begin{gather}
 H_{q}(\Phi_2 u)+\Psi_2 u+B_2(x^{-1}u h_q u)=0
,\label{2.8}
\end{gather}
where $\Phi_1$, $\Psi_1$, $B_1$, $\Phi_2$, $\Psi_2$, $B_2$ are
polynomials, $\Phi_1$, $\Phi_2$ monic and $\deg \Phi_i=t_{i}$,
$\deg \Psi_i=p_{i}$, $\deg B_i =r_{i}$, $d_{i}=\max(t_{i},r_{i})$, $s_{i}=\max(p_{i}-1,d_{i}-2)$
for $i \in \{1,2\}$. Let $\Phi=\gcd(\Phi_1, \Phi_2)$. Then, there
exist two polynomials $\Psi$ and $B$ such that
\begin{gather}
 H_{q}(\Phi u)+\Psi u+B\big(x^{-1}u h_q u\big)=0 ,\label{2.9}
\end{gather}
with
\begin{gather}
 s=\max (p-1, d-2)=
s_{1}-t_{1}+t=s_{2}-t_{2}+t,\label{2.10}
\end{gather}
where $t=\deg \Phi$, $ p=\deg \Psi$, $r=\deg B$ and $d=\max(t,r)$.
\end{proposition}

\begin{proof} With $\Phi=\gcd(\Phi_1,\Phi_2)$, there exist two co-prime polynomials
$\widetilde{\Phi}_1$, $\widetilde{\Phi}_2 $ such that
\begin{gather}
\Phi_1=\Phi\widetilde{\Phi}_1, \qquad
\Phi_2=\Phi \widetilde{\Phi}_2. \label{2.11}
\end{gather}
Taking into account \eqref{1.11} equations \eqref{2.7}, \eqref{2.8} become for $i \in
\{1,2\}$
\begin{gather}\label{2.12}
 \big(h_{q^{-1}}\widetilde{\Phi}_i\big)   H_{q}(\Phi
u)+\bigl\{\Psi_{i}+q^{-1} H_{q^{-1}}\widetilde{\Phi}_{i} \bigr\}
u+B_{i} \big(x^{-1}u h_q u\big)=0.
\end{gather}
The operation $(h_{q^{-1}}\widetilde{\Phi}_2)\times
(\ref{2.12}_{i=1})-(h_{q^{-1}}\widetilde{\Phi}_1)\times (\ref{2.12}_{i=2})$ gives
\begin{gather*}
\bigl\{\big(h_{q^{-1}} \widetilde{\Phi}_2\big)\bigl(\Psi_1 + q^{-1} \Phi
\big(H_{q^{-1}} \widetilde{\Phi}_1\big) \bigr)-\big(h_{q^{-1}}
\widetilde{\Phi}_1\big)\bigl(\Psi_2 + q^{-1} \Phi \big(H_{q^{-1}}
\widetilde{\Phi}_2\big)\bigr)\bigr\}u
\\
\qquad
 +\bigl\{\big(h_{q^{-1}} \widetilde{\Phi}_2\big) B_1-\big(h_{q^{-1}}
\widetilde{\Phi}_1\big) B_2\bigr\}\big(x^{-1}uh_q u\big)=0.
\end{gather*}
From the fact that $u$ is a strict $q$-Laguerre--Hahn form and by
virtue of Lemma~\ref{lemma3}  we get
\begin{gather*}
\big(h_{q^{-1}} \widetilde{\Phi}_1\big)\bigl(\Psi_2 +q^{-1} \Phi \big(H_{q^{-1}}
\widetilde{\Phi}_2\big)\bigr)=\big(h_{q^{-1}} \widetilde{\Phi}_2\big)\bigl(
\Psi_1 +q^{-1}\Phi \big(H_{q^{-1}} \widetilde{\Phi}_1\big)\bigr)   , \\
\big(h_{q^{-1}}\widetilde{\Phi}_1\big)B_2=\big(h_{q^{-1}}\widetilde{\Phi}_2\big)B_1.
\end{gather*}
Thus, there exist two polynomials $\Psi$ and $B$ such that
\begin{gather}
\Psi_1+ q^{-1}\Phi
\big(H_{q^{-1}}\widetilde{\Phi}_1\big)=\big(h_{q^{-1}}\widetilde{\Phi}_1\big)\Psi ,\qquad
 \Psi_2+ q^{-1} \Phi \big(H_{q^{-1}}\widetilde{\Phi}_2\big)=\big(h_{q^{-1}} \widetilde{\Phi}_2\big)\Psi ,\nonumber\\
 B_1=(h_{q^{-1}} \widetilde{\Phi}_1) B ,\qquad
 B_2=(h_{q^{-1}}\widetilde{\Phi}_2) B .
\label{2.13}
\end{gather}
Then, formulas \eqref{2.7}, \eqref{2.8} become
\begin{gather}\label{2.14}
\big(h_{q^{-1}} \widetilde{\Phi}_i\big)\bigl\{H_{q}(\Phi u)+\Psi
u+B\big(x^{-1}u h_q u\big)\bigr\}=0 ,\qquad i\in
\{1,2\}.
\end{gather}
But the polynomials $h_{q^{-1}}\widetilde{\Phi}_1$ and
$h_{q^{-1}}\widetilde{\Phi}_2$ are also co-prime. Using the Bezout
identity, there exist two polynomials $A_{1}$ and $A_{2}$ such that
\[
A_{1} \big(h_{q^{-1}}\widetilde{\Phi}_1\big)+A_{2}
\big(h_{q^{-1}}\widetilde{\Phi}_2\big)=1.
\]
Consequently, the operation $A_{1}\times (\ref{2.14}_{i=1})+A_{2}\times
(\ref{2.14}_{i=2})$ leads to \eqref{2.9}. With~\eqref{2.11} and~\eqref{2.13} it is easy to
prove~\eqref{2.10}.
\end{proof}

\begin{proposition}\label{proposition3} For any $q$-Laguerre--Hahn form $u$, the
triplet $(\Phi , \Psi, B)$ $(\Phi$ monic$)$ which realizes the minimum
of $\hbar(u)$ is unique.
\end{proposition}

\begin{proof} If $s_{1}=s_{2}$ in \eqref{2.9}, \eqref{2.10} and $s_{1}=s_{2}=s=\min
\hbar(u)$, then $t_{1}=t=t_{2}$. Consequently,
$\Phi_{1}=\Phi=\Phi_{2}$,   $B_{1}=B=B_{2}$ and
$\Psi_{1}=\Psi=\Psi_{2}$.
\end{proof}

Then, it's necessary to give a criterion which allows us to simplify
the class. For this, let us recall the following lemma:
\begin{lemma}\label{lemma5}
Consider $u$ a regular form, $\Phi$, $\Psi $ and $B$ three
polynomials, $\Phi$ monic. For any zero~$c$ of~$\Phi$, denoting
\begin{gather}
\Phi (x)=(x-c)\Phi_{c}(x),\nonumber
 \\
 q\Psi (x)+ \Phi_{c}(x )=(x-cq)\Psi_{cq}(x )+r_{cq},\label{2.15}\\
 qB(x)=(x-cq)B_{cq}(x )+b_{cq} .\nonumber
\end{gather}
The following statements are equivalent:
\begin{gather}
H_{q}(\Phi u)+\Psi u+B\big(x^{-1}u h_q u\big)=0,\nonumber
\\
H_{q}(\Phi_{c} u)+\Psi_{cq} u+B_{cq}\big(x^{-1}u h_q
u\big)+r_{cq}(x-cq)^{-1}u+b_{cq}(x-cq)^{-1}\big(x^{-1}u h_q u\big)
\nonumber\\
\qquad{} -\bigl\{\langle u,\Psi_{cq} \rangle + \langle x^{-1}u
h_q u,B_{cq} \rangle \bigr\}\delta_{cq}=0. \label{2.16}
\end{gather}
\end{lemma}

\begin{proof}The proof is obtained straightforwardly by using the relations in
\eqref{1.2} and in~\eqref{2.1}.
\end{proof}

\begin{proposition}\label{proposition4} A regular form $u$ $q$-Laguerre--Hahn
satisfying \eqref{1.17} is of class $s$ if and only if
\begin{gather}
\prod_{c \in \mathcal{Z}_{\Phi}} \Big\{ |q (h_{q}\Psi)(c)+
(H_{q}\Phi)(c) |+ |q (h_{q}B)(c) |\nonumber\\
\qquad{} +\big|\langle
u,q (\theta_{cq}\Psi)+(\theta_{cq}\circ \theta_{c}\Phi)+q
 (h_{q}u(\theta_{0}\circ
\theta_{cq}B) )\rangle\big|\Big\}>0, \label{2.17}
\end{gather}
where $\mathcal{Z}_{\Phi}$ is the set of roots of $\Phi$.
\end{proposition}

\begin{proof} Let $c$ be a root of $\Phi$: $\Phi(x)=(x-c)
\Phi_{c}(x)$. On account of \eqref{2.15} we have
\begin{gather*}
r_{cq}=q \Psi(cq)+\Phi_{c}(cq)=q (h_{q}\Psi)(c)+ (H_{q}\Phi)(c),\qquad
  b_{cq}=q B(cq)=q (h_{q}B)(c),
\\
\Psi_{cq}(x)=q (\theta_{cq}\Psi)(x)+(\theta_{cq}\Phi_{c})(x)=q
(\theta_{cq}\Psi)(x)+(\theta_{cq}\circ\theta_{c}\Phi)(x),
\\
 B_{cq}(x)=q (\theta_{cq}B)(x).
\end{gather*}
Therefore,
\begin{gather*}
   \langle u,\Psi_{cq} \rangle + \langle x^{-1}u
h_q u,B_{cq} \rangle  =  \langle
u,q(\theta_{cq}\Psi)+(\theta_{cq}\circ\theta_{c}\Phi)\rangle
+\langle u h_{q}u,q \theta_{0}\circ\theta_{cq}B\rangle \\
\phantom{\langle u,\Psi_{cq} \rangle + \langle x^{-1}u h_q u,B_{cq} \rangle}{}
=  \langle u,q(\theta_{cq}\Psi)+(\theta_{cq}\circ\theta_{c}\Phi)\rangle+
\langle u,q (h_{q}u(\theta_{0}\circ\theta_{cq}B))\rangle
 \\
\phantom{\langle u,\Psi_{cq} \rangle + \langle x^{-1}u h_q u,B_{cq} \rangle}{}
   =\langle u,q(\theta_{cq}\Psi)+(\theta_{cq}\circ\theta_{c}\Phi)+
q (h_{q}u(\theta_{0}\circ\theta_{cq}B))\rangle.
\end{gather*}
The condition \eqref{2.17} is necessary. Let us suppose that $c$ fulf\/ils
the conditions
\begin{gather*}
r_{cq}=0, \qquad b_{cq}=0, \qquad \langle
u,q(\theta_{cq}\Psi)+(\theta_{cq}\circ\theta_{c}\Phi)+ q
(h_{q}u(\theta_{0}\circ\theta_{cq}B))\rangle=0.
\end{gather*}
Then on account of Lemma~\ref{lemma5} \eqref{2.16} becomes
\[
H_{q}(\Phi_{c} u)+\Psi_{cq} u+B_{cq}\big(x^{-1}u h_q u\big)=0
\]
with $s_{c}=\max(\max(\deg \Phi_{c},\deg B_{cq})-2,\deg
\Psi_{c}-1) <s$, what contradicts with $s:=\min \hbar(u)$.

The condition \eqref{2.17} is suf\/f\/icient. Let us suppose $u$ to be of
class $\widetilde{s}<s$. There exist three polynomials
$\widetilde{\Phi}$ (monic) $\deg \widetilde{\Phi}=\widetilde{t}$,
$\widetilde{\Psi}$, $\deg \widetilde{\Phi}=\widetilde{p}$,
$\widetilde{B}$, $\deg\widetilde{B}=\widetilde{r}$ such that
\[
H_{q}(\widetilde{\Phi} u)+\widetilde{\Psi} u+\widetilde{B}(x^{-1}u
h_q u)=0
\]
with $\widetilde{s}=\max(\widetilde{d}-2,\widetilde{p}-1)$ where
$\widetilde{d}:=\max(\widetilde{t},\widetilde{r})$. By Proposition~\ref{proposition2}, it exists a polynomial $\chi$ such that
\begin{gather*}
\Phi=\chi \, \widetilde{\Phi}, \qquad
\Psi=(h_{q^{-1}}\chi)\widetilde{\Psi}-q^{-1}
(H_{q^{-1}}\chi)\widetilde{\Phi} , \qquad
B=(h_{q^{-1}}\chi)\widetilde{B}.
\end{gather*}
Since $\widetilde{s}<s$ hence $\deg \chi \geq 1$. Let $c$ be a zero
of $\chi: \chi(x)=(x-c) \chi_{c}(x)$. On account of \eqref{1.10} we
have
\[
q \Psi(x)+\Phi_{c}(x)=(x-cq) \bigl\{(h_{q^{-1}}\chi_{c})(x)
\widetilde{\Psi}(x)-q^{-1}
(H_{q^{-1}}\chi_{c})(x)\widetilde{\Phi}(x) \bigr\}.
\]
Thus $r_{cq}=0$ and  $b_{cq}=0$. Moreover, with \eqref{1.8} we have
\begin{gather*}
\bigl\langle u,q (\theta_{cq}\Psi)+(\theta_{cq}\circ
\theta_{c}\Phi)+ q \bigl(h_{q}u(\theta_{0}\circ
\theta_{cq}B)\bigr)\bigr\rangle \\
\qquad{}  = \bigl\langle
u,(h_{q^{-1}}\chi_{c}) \widetilde{\Psi}-q^{-1}
(H_{q^{-1}}\chi_{c})\widetilde{\Phi}+(h_{q}u)\theta_{0}((h_{q^{-1}}\chi_{c})
\widetilde{B})\bigr\rangle\\
\qquad{} = \bigl\langle u,(h_{q^{-1}}\chi_{c})
\widetilde{\Psi}- (H_{q}\circ
h_{q^{-1}}\chi_{c})\widetilde{\Phi}+(h_{q}u)\theta_{0}((h_{q^{-1}}\chi_{c})
\widetilde{B})\bigr\rangle\\
\qquad{} = \langle \widetilde{\Psi}u,h_{q^{-1}}\chi_{c}\rangle+\langle H_{q}(\widetilde{\Phi}u),h_{q^{-1}}\chi_{c}\rangle
+\langle \widetilde{B}\big(x^{-1}u h_{q}u\big),h_{q^{-1}}\chi_{c}\rangle  \\
\qquad{} = \bigl\langle H_{q}(\widetilde{\Phi}u)+\widetilde{\Psi}u+\widetilde{B}\big(x^{-1}u h_{q}u\big),h_{q^{-1}}\chi_{c}\bigr\rangle
 = 0.
\end{gather*}
This is contradictory with \eqref{2.17}. Consequently, $\widetilde{s}=s$,
$\widetilde{\Phi}=\Phi$, $\widetilde{\Psi}=\Psi$ and
$\widetilde{B}=B$.
\end{proof}

\begin{remark}\label{remark2} When $q\longrightarrow 1$ we recover again the
criterion which allows us to simplify a usual Laguerre--Hahn form~\cite{6}.
\end{remark}

\begin{remark}\label{remark3} When $B=0$ and $s=0$, the form $u$ is usually called
$q$-classical~\cite{16}. When $B=0$ and $s=1$, the symmetrical
$q$-semiclassical orthogonal $q$-polynomials of class one are
exhaustively described in~\cite{14}.
\end{remark}

\begin{proposition}\label{proposition5}
Let $u$ be a symmetrical $q$-Laguerre--Hahn form of class $s$
satisfying~\eqref{1.17}. The following statements hold
 \begin{itemize}\itemsep=0pt
\item[$(i)$]If  $s$ is odd, then the polynomials~$\Phi$ and $B$ are odd and~$\Psi$ is even.
\item[$(ii)$]If  $s$ is even, then the polynomials~$\Phi$ and~$B$ are even and~$\Psi$ is odd.
 \end{itemize}
\end{proposition}

\begin{proof}
Writing
\begin{gather*}
\Phi(x)=\Phi^{\rm e}\big(x^2\big)+x\Phi^{\rm o}\big(x^2\big) ,\qquad\!\! \Psi(x)=\Psi^{\rm e}\big(x^2\big)+x\Psi^{\rm o}\big(x^2\big)
 ,\qquad\!\! B(x)=B^{\rm e}\big(x^2\big)+xB^{\rm o}\big(x^2\big) ,
\end{gather*}
then \eqref{1.17} becomes
\begin{gather*}
H_{q}\big(\Phi^{\rm e}\big(x^2\big) u\big)+x\Psi^{\rm o}\big(x^2\big) u+B^{\rm e}\big(x^2\big)\big(x^{-1}u h_q
u\big)\\
\qquad{} +H_{q}\big(x\Phi^{\rm o}\big(x^2\big) u\big)+\Psi^{\rm e}\big(x^2\big) u+xB^{\rm o}\big(x^2\big)\big(x^{-1}u h_q u\big)=0.
\end{gather*}
Denoting
\begin{gather}
w^{\rm e}=H_{q}\big(\Phi^{\rm e}\big(x^2\big) u\big)+x\Psi^{\rm o}\big(x^2\big)
u+B^{\rm e}\big(x^2\big)\big(x^{-1}u h_q u\big),\nonumber\\
w^{\rm o}=H_{q}\big(x\Phi^{\rm o}\big(x^2\big) u\big)+\Psi^{\rm e}\big(x^2\big) u+xB^{\rm o}\big(x^2\big)\big(x^{-1}u h_q u\big).
\label{2.18}
\end{gather}
Then,
\begin{gather}
w^{\rm o}+w^{\rm e}=0. \label{2.19}
\end{gather}
From \eqref{2.19} we get
\begin{gather}
(w^{\rm o})_n=-(w^{\rm e})_n , \qquad  n\geq0. \label{2.20}
\end{gather}
From def\/initions in \eqref{2.18} and \eqref{2.20} we can write for    $n\geq0$
\begin{gather}
(w^{\rm e})_{2n}=\langle u, x^{2n+1}\Psi^{\rm o}\big(x^2\big)-[2n]_q
x^{2n-1}\Phi^{\rm e}\big(x^2\big)\rangle +\langle uh_qu,x^{2n-1}
B^{\rm e}\big(x^2\big)\rangle , \nonumber\\
(w^{\rm o})_{2n+1}=\langle u,x^{2n+1}\Psi^{\rm e}\big(x^2\big)-[2n+1]_q
x^{2n+1}\Phi^{\rm o}\big(x^2\big)\rangle +\langle uh_q u,x^{2n+1}
B^{\rm o}\big(x^2\big)\rangle.
\label{2.21}
\end{gather}
Now, with the fact that~$u$ is a symmetrical form then~$uh_{q}u$ is
also a symmetrical form. Indeed,
\begin{gather*}
   (uh_{q}u)_{2n+1}=\sum_{k=0}^{2n+1} (h_{q}u)_{k} (u)_{2n+1-k} =  \sum_{k=0}^{2n+1} q^{k}(u)_{k} (u)_{2n+1-k}\\
\phantom{(uh_{q}u)_{2n+1}}{} =  \sum_{k=0}^{n} q^{2k}(u)_{2k} (u)_{2(n-k)+1}+\sum_{k=0}^{n}
q^{2k+1}(u)_{2k+1}(u)_{2(n-k)} = 0, \qquad n\geq 0.
\end{gather*}
Thus \eqref{2.21} gives
\begin{gather}
(w^{\rm o})_{2n+1}=0=(w^{\rm e})_{2n}   , \qquad n\geq0.\label{2.22}
\end{gather}
On account of \eqref{2.19} and \eqref{2.22} we deduce   $w^{\rm o}=w^{\rm e}=0$.
Consequently $u$ satisf\/ies two $q$-dif\/ference equations
\begin{gather}
H_{q}\big(\Phi^{\rm e}\big(x^2\big) u\big)+x\Psi^{\rm o}\big(x^2\big) u+B^{\rm e}\big(x^2\big)\big(x^{-1}u h_q
u\big)=0,\label{2.23}
\end{gather}
and
\begin{gather}
H_{q}\big(x\Phi^{\rm o}\big(x^2\big) u\big)+\Psi^{\rm e}\big(x^2\big) u+xB^{\rm o}\big(x^2\big)\big(x^{-1}u h_q
u\big)=0.\label{2.24}
\end{gather}

$(i)$ If  $s=2k+1$,  with   $s=\max(d-2,p-1)$  we get
 $d \leq 2k+3$, $p\leq 2k+2$  then   $\deg (x\Psi^{\rm o}(x^2))\leq
2k+1$, $\deg (\Phi^{\rm e}(x^2))\leq 2k+2$  and  $\deg (B^{\rm e}(x^2))\leq
2k+2$. So, in accordance with \eqref{2.23}, we obtain the contradiction
 $s=2k+1\leq 2k$. Necessary  $\Phi^{\rm e}=B^{\rm e}=\Psi^{\rm o}=0$.

$(ii)$ If  $s=2k$,  with   $s=\max(d-2,p-1)$  we get  $d
\leq 2k+2$, $p\leq 2k+1$  then   $\deg (\Psi^{\rm e}(x^2))\leq
2k$, $\deg (x\Phi^{\rm o}(x^2))\leq 2k+1$  and  $\deg (x B^{\rm o}(x^2))\leq
2k+1$. So, in accordance with \eqref{2.24}, we obtain the contradiction
 $s=2k\leq 2k-1$. Necessary   $\Phi^{\rm o}=B^{\rm o}=\Psi^{\rm e}=0$. Hence the
desired result.
\end{proof}

\section[Different characterizations of $q$-Laguerre-Hahn forms]{Dif\/ferent characterizations of $\boldsymbol{q}$-Laguerre--Hahn forms}

One of the most important characterizations of the $q$-Laguerre--Hahn
forms is given in terms of a non homogeneous second order
$q$-dif\/ference equation so called $q$-Riccati equation fulf\/illed by
its formal Stieltjes series. See also \cite{6,8,10,15} for the usual case
and~\cite{13} for the discrete one.
\begin{proposition}\label{proposition6}
Let $u $ be a regular form. The following statement are equivalents:
\begin{itemize}\itemsep=0pt
\item[$(a)$] $u$ belongs to the $q$-Laguerre--Hahn class, satisfying \eqref{1.17}.
\item[$(b)$] The Stieljes formal series $S(u)$ satisfies the
 $q$-Riccati equation
\begin{gather}
(h_{q^{-1} }\Phi)(z) H_{q ^{-1}}(S( u))(z)= B(z) S(u)(z) (h_{q^{-1}}
S(u))(z)+C(z) S(u)(z) + D(z), \label{3.1}
\end{gather}
where  $\Phi$ and $B$ are polynomials defined in \eqref{1.17} and
\begin{gather}
C(z)=-(H_{q ^{-1}}\Phi)(z) - q \Psi(z),\nonumber\\
D(z)=-\bigl\{H_{q ^{-1}}(u\theta_0\Phi)(z)+q (u\theta_0\Psi)(z) +
q(u h_q u)\big(\theta_0^2 B\big)(z)\bigr\}.
\label{3.2}
\end{gather}
\end{itemize}
\end{proposition}

\begin{proof}
$(a)$ $\Rightarrow$ $(b)$. Suppose that $(a)$
is satisf\/ied, then there exist three polynomials   $\Phi$
(monic),   $\Psi$  and  $B$  such that  $H_q (\Phi u)+\Psi u +
B(x^{-1}uh_q u)=0$. From \eqref{1.11} the above $q$-dif\/ference equation
becomes
\[
(h_{q^{-1}}\Phi) (H_q u)+\bigl\{\Psi+q^{-1}(H_{q^{-1}}\Phi)\bigr \}
u + B\big(x^{-1}uh_q u\big)=0.
\]
From def\/inition of $S(u)$ and the linearity of $S$ we obtain
\begin{gather}
S\bigl((h_{q^{-1}}\Phi) (H_q u)\bigr)(z)+S(\Psi
u)(z)+q^{-1}S\bigl((H_{q^{-1}}\Phi) u \bigr)(z)+ S\bigl(B(x^{-1}uh_q
u)\bigr)(z)=0.\label{3.3}
\end{gather}
Moreover,
\begin{gather*}
  S(\Psi u)(z) \overset{\rm by \ \eqref{1.13}}{=} \Psi(z) S(u)(z)+(u\theta_{0}\Psi)(z),\\
q^{-1} S((H_{q^{-1}}\Phi) u)(z) \overset{\rm by \ \eqref{1.13}}{=}  q^{-1}(H_{q^{-1}}\Phi)(z)
S(u)(z)+q^{-1}(u\theta_{0}(H_{q^{-1}}\Phi))(z),
\\
S\bigl((h_{q^{-1}}\Phi) (H_q u)\bigr)(z)  \overset{\rm by \ \eqref{1.13}}{=}  (h_{q^{-1}}\Phi)(z)
S(H_q u)(z)+
\bigl((H_q u)\theta_{0}(h_{q^{-1}}\Phi)\bigr)(z)\\
\phantom{S\bigl((h_{q^{-1}}\Phi) (H_q u)\bigr)(z)}{}
\overset{\rm by \ \eqref{1.16}}{=}  q^{-1} (h_{q^{-1}}\Phi)(z) H_{q ^{-1}}(S( u))(z)+
\bigl((H_q u)\theta_{0}(h_{q^{-1}}\Phi)\bigr)(z),\\
S\bigl(B(x^{-1}uh_q u)\bigr)(z) \overset{\rm by \ \eqref{1.13}}{=}  B(z)S\bigl(x^{-1}uh_q
u)\bigr)(z)+\bigl((x^{-1}uh_q u)\theta_{0}B\bigr)(z) \\
\phantom{S\bigl(B(x^{-1}uh_q u)\bigr)(z)}{}
    \overset{\rm by \ \eqref{1.15}}{=}  z^{-1} B(z) S\bigl(uh_q
u)\bigr)(z)+\bigl((uh_q u)\theta_{0}^{2}B\bigr)(z) \\
\phantom{S\bigl(B(x^{-1}uh_q u)\bigr)(z)}{}
    \overset{\rm by \ \eqref{1.14}}{=}
      -B(z)S(u)(z) S(h_q u)(z)+ \bigl((uh_q u)\theta_{0}^{2}B\bigr)(z)\\
\phantom{S\bigl(B(x^{-1}uh_q u)\bigr)(z)}{}
    \overset{\rm by \ \eqref{1.16}}{=}
      -q^{-1}B(z) S(u)(z) (h_{q^{-1}}S(u))(z)+ \bigl((uh_q u)\theta_{0}^{2}B\bigr)(z), 
\end{gather*}
and
\[
(u\theta_{0}(H_{q^{-1}}\Phi))(z)+q \bigl((H_q
u)\theta_{0}(h_{q^{-1}}\Phi)\bigr)(z) \overset{\rm by \ \eqref{1.12}}{=} H_{q^{-1}}(u \theta_{0}
\Phi)(z).
\]
\eqref{3.3} becomes
\begin{gather*}
(h_{q^{-1} }\Phi)(z) H_{q ^{-1}}(S( u))(z)= B(z)S(u)(z) (h_{q^{-1} }
S(u))(z)-(H_{q ^{-1}}\Phi + q\Psi)(z) S(u)(z)\\
\phantom{(h_{q^{-1} }\Phi)(z) H_{q ^{-1}}(S( u))(z)=}{} - \bigl\{H_{q^{-1}}(u \theta_{0} \Phi)+q u\theta_0\Psi
+ q (u h_q u)\theta_0^2 B \bigr\}(z).
\end{gather*}
The previous relation gives \eqref{3.1} with
\eqref{3.2}.

$(b)$ $\Rightarrow$ $(a)$. Let $u \in \mathcal{P}'$
regular with its formal Stieltjes series $S(u)$ satisfying~\eqref{3.1}.
Likewise as in the previous implication, formula~\eqref{3.1} leads to
\begin{gather*}
S\bigl\{H_q (\Phi u)-q ^{-1}(C+H_{q ^{-1}}\Phi)u+B\big(x^{-1}uh_q
u\big)\bigr\}
\\
\qquad {}=
  q^{-1}D-q^{-1} u\theta_0 C + \big((u h_q u)\theta_0^2
B\big)+((H_{q}u)\theta_{0}(h_{q ^{-1}}\Phi)),
\end{gather*}
which implies
\begin{gather*}
  S\bigl\{H_q (\Phi u)-q^{-1}(C+H_{q ^{-1}}\Phi)u+B\big(x^{-1}uh_q
u\big)\bigr\}=0 ,\\
 D(z)=(u\theta_0 C)(z) -q \big((u h_q u)\big(\theta_0^2 B\big)\big)(z)-q
((H_{q}u)\theta_{0}(h_{q ^{-1}}\Phi))(z).
\end{gather*}
According to \eqref{3.2} and \eqref{1.12} we deduce that
\[
H_{q}(\Phi u)+\Psi u+B\big(x^{-1}u h_q u\big)=0 ,
\]
with
\begin{gather}
\Psi=-q^{-1}(C+H_{q^{-1}}\Phi).\label{3.4}
\end{gather}
\end{proof}

We are going to give the criterion which allows us to simplify the
class of $q$-Laguerre--Hahn form in terms of the coef\/f\/icients
corresponding to the previous characterization.
\begin{proposition}\label{proposition7} A regular form $u$ $q$-Laguerre--Hahn
satisfying \eqref{3.1} is of class $s$ if and only if
\begin{gather}
\prod\limits_{c\in Z_\Phi} \bigl\{|B(cq)|+|C(cq)|+ |D(cq)|\bigr\}>
0, \label{3.5}
\end{gather}
where $Z_\Phi$ is the set of roots of $\Phi$ with
\begin{gather}
s=\max\bigl(\deg B-2, \deg C-1, \deg D\bigr). \label{3.6}
\end{gather}
\end{proposition}

\begin{proof} By comparing \eqref{2.17} and \eqref{3.5}, it is enough to prove the following equalities
\begin{gather*}
|C(cq)|=\bigl|q (h_{q}\Psi)(c)+ (H_{q}\Phi)(c)\bigr| ,
\\
|D(cq)|=\bigl|\bigl\langle u,q (\theta_{cq}\Psi)+(\theta_{cq}\circ
\theta_{c}\Phi)+q \bigl(h_{q}u(\theta_{0}\circ
\theta_{cq}B)\bigr)\bigr\rangle\bigr|.
\end{gather*}
Indeed, on account of \eqref{3.2}, the def\/inition of the polynomial $u f$,
the def\/inition of the product form $uv$ and~\eqref{1.8} we have
\[
C(cq)=-(H_{q ^{-1}}\Phi)(cq) - q \Psi(cq)=-(H_{q}\Phi)(c)- q
(h_{q}\Psi)(c),
\]
and
\begin{gather*}
  D(cq) = -\bigl\{H_{q ^{-1}}(u\theta_0\Phi)(cq)+q (u\theta_0\Psi)(cq) +
q(u h_q u)({\theta_0}^2 B)(cq)\bigr\} \\
\phantom{D(cq)}{}=  -\bigl\{H_{q}(u\theta_0\Phi)(c)+\langle u,q \theta_{cq}\Psi\rangle +
\langle u h_q u,q \theta_0\circ\theta_{cq} B\rangle\bigr\}\\
\phantom{D(cq)}{}= -\bigl\{H_{q}(u\theta_0\Phi)(c)+\bigl\langle u,q
\theta_{cq}\Psi+q \bigl(h_{q}u(\theta_{0}\circ
\theta_{cq}B)\bigr)\bigr\rangle \bigr\}.
\end{gather*}
Moreover,
\begin{gather*}
  H_{q}(u\theta_0\Phi)(c) \overset{\rm by \ \eqref{1.6}}{=} \frac{(u\theta_0\Phi)(cq)-(u\theta_0\Phi)(c)}{(q-1)c}
  = \bigl\langle u,
\frac{\theta_{cq}\Phi-\theta_{c}\Phi}{cq-c}\bigr\rangle =
\bigl\langle u, \theta_{cq}\circ \theta_{c} \Phi\bigr\rangle.
\end{gather*}
Thus \eqref{2.17} is equivalent to \eqref{3.5}. To prove \eqref{3.6}, according to the
def\/inition of the class we may write
\begin{gather}
s=\max\bigl(\deg B-2, \deg \Phi-2, \deg \Psi-1\bigr). \label{3.7}
\end{gather}

$\bullet$ If $\deg \Psi\neq \max\bigl(\deg B-1, \deg \Phi-1\bigr)$,
on account of \eqref{3.2} and \eqref{3.7} we get the following implications
\begin{gather*}
\deg B\leq \deg \Phi \Rightarrow \left\{
                                   \begin{array}{l}
                                     \deg C=s+1,\\
                                     \deg D\leq s
                                   \end{array}
                                 \right.
\Rightarrow \max\bigl(\deg B-2, \deg C-1, \deg D\bigr)=s,
\\
\deg B> \deg \Phi \Rightarrow \left\{
                                   \begin{array}{l}
                                     \deg C\leq s+1,\\
                                     \deg D= s
                                   \end{array}
                                 \right.
\Rightarrow \max\bigl(\deg B-2, \deg C-1, \deg D\bigr)=s.
\end{gather*}

$\bullet$ If $\deg \Psi= \max\bigl(\deg B-1, \deg \Phi-1\bigr)$ and
$\deg B> \deg \Phi$ then $s+1=\deg \Psi=\deg B-1> \deg \Phi-1$.
Consequently, $\max\bigl(\deg B-2, \deg C-1, \deg D\bigr)=s$.

$\bullet$ If $\deg \Psi= \max\bigl(\deg B-1, \deg \Phi-1\bigr)$ and
$\deg B=\deg \Phi$ then $\deg \Psi=\deg B-1=\deg \Phi-1$ which
implies $\deg B-2=s$, $\deg C-1\leq s$, $\deg D\leq s$. Therefore
$\max\bigl(\deg B-2, \deg C-1, \deg D\bigr)=s$.

$\bullet$ If $\deg \Psi= \max\bigl(\deg B-1, \deg \Phi-1\bigr)$ and
$\deg B < \deg \Phi$ then $\deg \Psi=\deg \Phi-1$ and $s=\deg
\Psi-1$. Writing $\Phi(x)= x^{p+1}+\textrm{lower order terms},
\Psi(x)=a_{p}x^{p}+\cdots +a_{0}$, by virtue of \eqref{3.2} and \eqref{1.3}, it is
worth noting that $C(z)=-\bigl([p+1]_{q^{-1}}+q a_{p}\bigr)z^{p-1}+
\textrm{lower order terms}$ and $D(z)=-\bigl([p]_{q^{-1}}+q
a_{p}\bigr)z^{p-1}+ \textrm{lower order terms}$ with
$[p+1]_{q^{-1}}\neq [p]_{q^{-1}}$ assuming either $\deg C=s$ or
$\deg D=s$. Thus, $ \max\bigl(\deg B-2, \deg
C-1, \deg D\bigr)=s.$

Hence the desired result \eqref{3.6}.
\end{proof}

\smallskip

An other important characterization of the $q$-Laguerre--Hahn forms
is the structure relation. See also \cite{6,15} for the usual case and~\cite{13} for the discrete one.

\begin{proposition}\label{proposition8} Let $u$ be a regular form and $\{P_{n}\}_{n\geq
0}$ be its MOPS. The following statements are equivalent:
\begin{itemize}\itemsep=0pt
\item[$(i)$]
 $u$ is a $q$-Laguerre--Hahn form satisfying \eqref{1.17}.
\item[$(ii)$] There exist an integer $s\geq 0$, two polynomials $\Phi$
$($monic$)$, $B$ with $t=\deg\Phi\leq s+2$,   $r=\deg B\leq s+2$ and a
sequence of complex numbers $\{\lambda_{n,\nu}\}_{n,\nu \geq 0}$
such that
\begin{gather}
\Phi(x) (H_{q}P_{n+1})(x)-h_{q}(B P_{n}^{(1)})(x)=
\sum_{\nu=n-s}^{n+d}\lambda_{n,\nu}P_{\nu}(x), \qquad n> s ,\quad
\lambda_{n,n-s}\neq0,
  \label{3.8}
\end{gather}
where $d=\max(t,r)$ and $\big\{P_{n}^{(1)}\big\}_{n \geq 0}$ be the
associated sequence of the first kind for the sequence
$\{P_n\}_{n\geq0}$.
\end{itemize}
\end{proposition}

\begin{proof} $(i)$ $\Rightarrow$ $(ii)$. Beginning with the expression
$\Phi(x) (H_{q}P_{n+1})(x)-h_{q}\big(B P_{n}^{(1)}\big)(x)$ which is a
polynomial of degree at most $n+d$. Then, there exists a sequence of
complex numbers $\{\lambda_{n,\nu}\}_{n\geq 0, \, 0\leq \nu \leq
n+d}$ such that
\begin{gather}
\Phi(x) (H_{q}P_{n+1})(x)-(h_{q}B)(x)\big(h_{q}P_{n}^{(1)}\big)(x)=
\sum_{\nu=0}^{n+d}\lambda_{n,\nu}P_{\nu}(x), \qquad  n\geq0.
\label{3.9}
\end{gather}
Multiplying both sides of \eqref{3.9} by $P_{m}$, $0\leq m \leq n+d$ and
applying $u$ we get
\begin{gather}
\langle u, \Phi P_m (H_q P_{n+1})\rangle - \langle h_{q}u,
B(h_{q^{-1} }P_{m})(u\theta_0 P_{n+1})\rangle=\lambda_{n,m}\langle
u,P_{m}^{2}\rangle , \nonumber\\ n\geq 0   , \qquad 0\leq m \leq n+d.\label{3.10}
\end{gather}
On the other hand, applying $H_{q}(\Phi u)+\Psi
u+B(x^{-1}u h_q u)=0$ to $P_{n+1}(h_{q^{-1} }P_m)$, on account of
the def\/initions, \eqref{1.10} and \eqref{1.8} we obtain
\begin{gather*}
  0 = \langle H_{q}(\Phi u)+\Psi
u+B\big(x^{-1}u h_q u\big),P_{n+1}(h_{q^{-1} }P_m)\rangle \\
\phantom{0} =  \bigl\langle u,\Psi P_{n+1}(h_{q^{-1} }P_m)-\Phi H_{q}\bigl(P_{n+1}(h_{q^{-1} }P_m)\bigr)\bigr\rangle
+\bigl\langle  h_q u,u\theta_{0}(B P_{n+1}(h_{q^{-1}
}P_m))\bigr\rangle\\
\phantom{0}= \langle u,\big\{\Psi(h_{q^{-1} }P_m)-q^{-1}\Phi
(H_{q^{-1}}P_{m})\big\}P_{n+1}-\Phi P_{m}(H_{q}P_{n+1})\rangle\\
\phantom{0=}+\langle
h_q u,u\theta_{0}(B P_{n+1}(h_{q^{-1} }P_m))\rangle.
\end{gather*}
Thus, for $  n\geq 0$, $0\leq m \leq n+d$
\begin{gather}
\langle u,\Phi P_{m} (H_{q}P_{n+1})\rangle=\bigl\langle
u,\big\{\Psi(h_{q^{-1} }P_m)-q^{-1}\Phi
(H_{q^{-1}}P_{m})\big\}P_{n+1}\bigr\rangle\nonumber\\
\phantom{\langle u,\Phi P_{m} (H_{q}P_{n+1})\rangle=}{} +\langle h_q u,u\theta_{0}(B
P_{n+1}(h_{q^{-1} }P_m))\rangle.\label{3.11}
\end{gather}
Using \eqref{3.10}, \eqref{3.11} to eliminate $\langle u,\Phi P_{m}
(H_{q}P_{n+1})\rangle$ we get for $ n\geq 0$, $0\leq m \leq
n+d$
\begin{gather}
\bigl\langle u,\big\{\Psi(h_{q^{-1} }P_m)-q^{-1}\Phi
(H_{q^{-1}}P_{m})\big\}P_{n+1}\bigr\rangle \nonumber\\
\qquad{}
 +\bigl\langle h_q u,u\theta_{0}(B P_{n+1}(h_{q^{-1}
}P_m))-(h_{q^{-1}
}P_m)B(u\theta_{0}P_{n+1})\bigr\rangle=\lambda_{n,m}\langle
u,P_{m}^{2}\rangle.\label{3.12}
\end{gather}
Moreover, by virtue of \eqref{1.5} we have
$B(u\theta_{0}P_{n+1})=u\theta_{0}(BP_{n+1})$, $n>s$. Therefore,
taking into account \eqref{1.4} and def\/initions, \eqref{3.12}
yields for $n>s$, $0\leq m \leq n+d$
\begin{gather*}
\bigl\langle u,\big\{\Psi(h_{q^{-1} }P_m)-q^{-1}\Phi (H_{q^{-1}}P_{m})+B
((h_{q}u)\theta_{0}(h_{q^{-1}
}P_m))\big\}P_{n+1}\bigr\rangle=\lambda_{n,m}\langle u,P_{m}^{2}\rangle
\end{gather*}
with
\[
\deg \bigl\{\Psi (h_{q^{-1} }P_m)-q^{-1} \Phi (H_{q^{-1}
}P_m)+B((h_q u)\theta_0 ( h_{q^{-1} }P_m))\bigr\} \leq m+s+1.
\]
Consequently, the orthogonality of $\{P_n\}_{n\geq0}$ with respect
to $u$ gives
\[
\lambda_{n,m}=0, \qquad 0\leq m\leq n-s-1, \quad n\geq s+1,
\qquad \lambda_{n,n-s}\neq0.
\]
Hence the desired result \eqref{3.8}.

 $(ii)$ $\Rightarrow$ $(i)$. Let $v$ be the
form def\/ined by
\[
v:=H_q (\Phi u)+ B\big(x^{-1}uh_q u\big)+\left(\sum\limits_{i=0}^{s+1}a_i x^i\right)u
\]
with $a_i\in{\mathbb C}$, $ 0\leq i \leq s+1$. From
def\/initions and the hypothesis of $(ii)$ we may write successively
\begin{gather*}
   \langle v,P_{n+1}\rangle= \bigl\langle H_q (\Phi u)+ B\big(x^{-1}uh_q u\big),P_{n+1}\bigr\rangle+
\langle u,P_{n+1} \sum\limits_{i=0}^{s+1}a_i x^i\rangle \\
\phantom{\langle v,P_{n+1}\rangle}{} =  - \bigl\langle u,\Phi (H_{q}P_{n+1})-(h_q u)\theta_{0}(BP_{n+1})\bigr\rangle+
\langle u,P_{n+1} \sum\limits_{i=0}^{s+1}a_i x^i\rangle\\
\phantom{\langle v,P_{n+1}\rangle}{}= -\bigl\langle u, \sum_{\nu=n-s}^{n+d} \lambda_{n,\nu} P_{\nu}\bigr\rangle+
\langle u,P_{n+1} \sum\limits_{i=0}^{s+1}a_i x^i\rangle \\
\phantom{\langle v,P_{n+1}\rangle}{}= -\sum_{\nu=n-s}^{n+d} \lambda_{n,\nu} \langle
u,P_{\nu}\rangle+\sum\limits_{i=0}^{s+1}a_i \langle u,x^{i}
P_{n+1}\rangle, \qquad n>s.
\end{gather*}
From assumption of orthogonality of $\{P_{n}\}_{n \geq 0}$ with
respect to $u$ we get
\[
\langle v,P_{n}\rangle=0, \qquad n\geq s+2.
\]
In order to get  $\langle v,P_n\rangle=0,$ for any  $n\geq 0$, we
shall choose $ a_i$ with $i=0,1,\dots,s+1$, such that  $\langle
v,P_i\rangle=0,$ for $i=0,1,\dots,s+1$. These coef\/f\/icients $a_i$ are
determined in a unique way. Thus, we have deduced the existence of
polynomial $\Psi(x)= \sum\limits_{i=0}^{s+1}a_i x^i$ such that
$\langle v,P_n\rangle=0,$ for any  $n\geq 0$.  This leads to  $H_q
(\Phi u)+\Psi u +B(x^{-1}uh_q u)=0$ and the point $(i)$ is then
proved.
\end{proof}

\section{Applications}

\subsection[The co-recursive of a $q$-Laguerre-Hahn form]{The co-recursive of a $\boldsymbol{q}$-Laguerre--Hahn form}

Let $\mu$ be a complex number, $u$ a regular form and
$\{P_n\}_{n\geq0}$ be its corresponding MOPS sa\-tisfying~\eqref{1.1}. We
def\/ine the co-recursive $\big\{P_{n}^{[\mu]}\big\}_{n\geq0}$ of
$\{P_n\}_{n\geq0}$ as the family of monic polynomials satisfying the
following three-term recurrence relation \cite[Def\/inition~4.2]{20}
\begin{gather*}
   P_{0}^{[\mu]}(x)=1 ,\qquad  P_{1}^{[\mu]}(x)=x-\beta_{0}-\mu ,\nonumber\\
   P_{n+2}^{[\mu]}(x)=(x- \beta_{n+1})P_{n+1}^{[\mu]}(x)-\gamma_{n+1} P_{n}^{[\mu]}(x),\qquad  n \geq 0.
\end{gather*}
Denoting by $u^{[\mu]}$ its corresponding regular form. It is well
known that \cite[equation~(4.14)]{20}
\[
u^{[\mu]}=u   \bigl(\delta-\mu x ^{-1} u \bigr)^{-1}.
\]

\begin{proposition} \label{proposition9} If $u$ is a $q$-Laguerre--Hahn form
of class $s$, then $u^{[\mu]}$ is a $q$-Laguerre--Hahn form of the
same class~$s$.
\end{proposition}

\begin{proof} The relation linking $S(u)$ and $S(u^{[\mu]})$ is \cite[equation~(4.15)]{20}
$S(u^{[\mu]})= \frac{S(u)}{1+\mu S(u)}$ or equivalently
\begin{gather}
S(u)=\frac{S(u^{[\mu]})}{1-\mu S(u^{[\mu]})}.\label{4.2}
\end{gather}
From def\/initions and by virtue of \eqref{4.2} we have
\[
h_{q^{-1}}S(u)=\frac{h_{q^{-1}}S(u^{[\mu]})}{1-\mu
h_{q^{-1}}S(u^{[\mu]})}
\]
and
\begin{gather*}
   (H_{q^{-1}}S(u))(z) =  \frac{\frac{(h_{q^{-1}}S(u^{[\mu]}))(z)}{1-\mu
(h_{q^{-1}}S(u^{[\mu]}))(z)}-\frac{S(u^{[\mu]})(z)}{1-\mu S(u^{[\mu]})(z)}}{(q^{-1}-1)z} \\
 \phantom{(H_{q^{-1}}S(u))(z)}{} =  \frac{(H_{q^{-1}}S(u^{[\mu]}))(z)}{\bigl(1-\mu
(h_{q^{-1}}S(u^{[\mu]}))(z)\bigr)\bigl(1-\mu S(u^{[\mu]})(z)\bigr)}.
\end{gather*}
Replacing the above results in \eqref{3.1} the $q$-Riccati equation
becomes
\begin{gather*}
(h_{q ^{-1}}\Phi) \frac{H_{q ^{-1}}S(u^{[\mu]})}{\bigl(1-\mu
h_{q^{-1}}S(u^{[\mu]})\bigr)\bigl(1-\mu S(u^{[\mu]})\bigr)} \\
\\
\qquad{}=
 B   \frac{S(u^{[\mu]})}{1-\mu S(u^{[\mu]})}
\frac{h_{q^{-1}}S(u^{[\mu]})}{1-\mu h_{q^{-1}}S(u^{[\mu]})} +C
\frac{S(u^{[\mu]})}{1-\mu S(u^{[\mu]})} + D.
\end{gather*}
Equivalently
\begin{gather*}
(h_{q ^{-1}}\Phi)H_{q ^{-1}}S(u^{[\mu]})=
 B   S(u^{[\mu]}) h_{q^{-1}}S(u^{[\mu]})+ C   S(u^{[\mu]})
\bigl(1-\mu h_{q^{-1}}S(u^{[\mu]})\bigr)\\
\phantom{(h_{q ^{-1}}\Phi)H_{q ^{-1}}S(u^{[\mu]})=}{} +D \bigl(1-\mu
h_{q^{-1}}S(u^{[\mu]})\bigr)\bigl(1-\mu S(u^{[\mu]})\bigr).
\end{gather*}
 Therefore the $q$-Riccati equation satisf\/ied by
$S(u^{[\mu]})$
\begin{gather}
(h_{q ^{-1}}\Phi^{[\mu]} )H_{q ^{-1}}S(u^{[\mu]})=B^{[\mu]}
S(u^{[\mu]})h_{q ^{-1}}S(u^{[\mu]})+C^{[\mu]} S(u^{[\mu]}) +
D^{[\mu]} , \label{4.3}
\end{gather}
where
\begin{alignat}{3}
&  K \Phi^{[\mu]}(x)=\Phi(x)+\mu(1-q)x (h_q D)(x)  ,  \qquad &&
K B^{[\mu]}(x)=B(x)-\mu C(x)+\mu^2 D(x)  , & \nonumber \\
& K C^{[\mu]}(x )=C(x)-2\mu D(x)  ,\qquad &&
 K D^{[\mu]}(x)=D(x), & \label{4.4}
\end{alignat}
the non zero constant $K$ is chosen such that the polynomial
$\Phi^{[\mu]}$ is
monic.
$u^{[\mu]}$ is then a $q$-Laguerre--Hahn form.

On account of \eqref{3.2}, \eqref{3.4} and \eqref{4.4} we get
\begin{gather}
K \Psi^{[\mu]}=\Psi+\mu\big(q^{-1}D+h_{q}D\big).\label{4.5}
\end{gather}
As a consequence, the regular form $u^{[\mu]}$ fulf\/ils the following
$q$-dif\/ference equation
\begin{gather}
H_{q}\big(\Phi^{[\mu]} u^{[\mu]}\big)+\Psi^{[\mu]}
u^{[\mu]}+B^{[\mu]}\big(x^{-1}u^{[\mu]} h_q u^{[\mu]}\big)=0.\label{4.6}
\end{gather}
We suppose that the $q$-Riccati equation~\eqref{3.1} of $u$ is irreducible
of class~$s$. With respect to the class, we use the result~\eqref{3.5} of
Proposition~\ref{proposition7} and get for every zero $c$ of $\Phi^{[\mu]}$:
\begin{itemize}\itemsep=0pt
\item  If $D(cq)\neq0$, then $D^{[\mu]}(cq)=K^{-1}D(cq)\neq0$ and
equation \eqref{4.3} is not reducible.
\item We suppose that $D(cq)=0$. From the fact that
$\Phi^{[\mu]}(c)=0$, the f\/irst relation in~\eqref{4.4} leads to
$\Phi(c)=0$ and the third equality in
\eqref{4.4} gives $C^{[\mu]}(cq)=K^{-1} C(cq)$.
\end{itemize}
If $C(cq)\neq0$, then the equation \eqref{4.3} is still not reducible.
If $C(cq)=0=D(cq)$, then $B^{[\mu]}(cq)=K^{-1} B(cq)\neq 0$ since
$u$ is of
class~$s$.
We conclude that
\[
\big|B^{[\mu]}(cq)\big|+\big|C^{[\mu]}(cq)\big|+\big|D^{[\mu]}(cq)\big|> 0.
\]
Consequently, the class $s^{[\mu]}$ of $u^{[\mu]}$ is given by
$s^{[\mu]}=\max\bigl(\deg B^{[\mu]}-2,\deg C^{[\mu]}-1,\deg
D^{[\mu]}\bigr)$. Accordingly to the last equality in~\eqref{4.4} and
\eqref{3.6} we get $s^{[\mu]}=\max\bigl(\deg B^{[\mu]}-2,\deg
C^{[\mu]}-1$, $\deg D\bigr)$. A~discussion on the degree leads to
$s^{[\mu]}=s$.
\end{proof}

\begin{example}\label{example1} Let $u$ be a $q$-classical form satisfying
the $q$-analog of the distributional equation of Pearson type
\begin{gather}
H_{q}(\phi u)+\psi u=0, \label{4.7}
\end{gather}
where $\phi$ is a monic polynomial of degree at most two and $\psi$
a polynomial of degree one, the co-recursive $u^{[\mu]}$ of~$u$ is a
$q$-Laguerre--Hahn form of class zero. $u^{[\mu]}$ and the Stieltjes
function $S(u^{[\mu]})$ satisfy, respectively, the $q$-dif\/ference
equation~\eqref{4.6} and the $q$-Riccati equation~\eqref{4.3} where on account
of \eqref{4.4}, \eqref{4.5}
\begin{gather*}
 K \Phi^{[\mu]}(x)=\frac{\phi''(0)}{2}   x^{2}+\left\{\phi'(0)+\mu
(q-1)\left(\frac{\phi''(0)}{2}+q \psi'(0)\right)\right\} x+\phi(0)  ,
\\
 K \Psi^{[\mu]}(x)=\psi'(0) x+\psi(0)-\mu
\big(q^{-1}+1\big)\left(\frac{\phi''(0)}{2}+q \psi'(0)\right),
\\
 K B^{[\mu]}(x)=\mu \left\{\left(\big(q^{-1}+1\big) \frac{\phi''(0)}{2}+q
\psi'(0) \right)x+\phi'(0)+q \psi(0)-
\left(\frac{\phi''(0)}{2}+q \psi'(0)\right)\mu\right\}  , \\
 K C^{[\mu]}(x )= -\left(q\psi'(0)+(q^{-1}+1)
\frac{\phi''(0)}{2}\right)x-\phi'(0)-q \psi(0)+2\mu \left(\frac{\phi''(0)}{2}+q \psi'(0)\right) ,\\
 K D^{[\mu]}(x)=-\frac{\phi''(0)}{2}-q  \psi'(0) .
\end{gather*}
\end{example}

\subsection[The associated of a $q$-Laguerre-Hahn form]{The associated of a $\boldsymbol{q}$-Laguerre--Hahn form}

Let $u$ be a regular form and $\{P_n\}_{n\geq0}$ its corresponding
MOPS satisfying~\eqref{1.1}. The associated sequence of the f\/irst kind
$\big\{P_{n}^{(1)}\big\}_{n\geq0}$ of $\{P_n\}_{n\geq0}$ satisf\/ies the
following three-term recurrence relation~\cite{20}
\begin{gather*}
   P_{0}^{(1)}(x)=1 ,\qquad  P_{1}^{(1)}(x)=x-\beta_{1} , \nonumber\\
   P_{n+2}^{(1)}(x)=(x- \beta_{n+2})P_{n+1}^{(1)}(x)-\gamma_{n+2} P_{n}^{(1)}(x),\qquad  n \geq 0.
\end{gather*}
Denoting by $u^{(1)}$ its corresponding regular form.
\begin{proposition}\label{proposition10} If $u$ is a $q$-Laguerre--Hahn form
of class~$s$, then $u^{(1)}$ is a $q$-Laguerre--Hahn form of the same
class~$s$.
\end{proposition}

\begin{proof} We assume that the formal Stieltjes function $S(u)$ of $u$
satisf\/ies~\eqref{3.1}. The relationship between $S(u^{(1)})$ and $S(u)$ is
\cite[equation~(4.7)]{20}
\[
\gamma_{1}   S(u^{(1)})(z)=-\frac{1}{S(u)(z)}-(z-\beta_{0}).
\]
Consequently,
\begin{gather}
S(u)(z)=-\frac{1}{\gamma_{1}
S(u^{(1)})(z)+(z-\beta_{0})}.\label{4.9}
\end{gather}
From def\/initions and by virtue of~\eqref{4.9} we have
\[
h_{q^{-1}}(S(u))(z)=-\frac{1}{\gamma_{1}
h_{q^{-1}}(S(u^{(1)}))(z)+q^{-1}z-\beta_{0}}
\]
and
\[
H_{q^{-1}}(S(u))(z)=\frac{\gamma_{1}H_{q^{-1}}(S(u^{(1)}))(z)+1}{\bigl(\gamma_{1}
h_{q^{-1}}(S(u^{(1)}))(z)+q^{-1}z-\beta_{0}\bigr)\bigl(\gamma_{1}
S(u^{(1)})(z)+z-\beta_{0}\bigr)}.
\]
Substituting in~\eqref{3.1} the $q$-Riccati equation becomes
\begin{gather*}
(h_{q^{-1}}\Phi)(z)
\frac{\gamma_{1}H_{q^{-1}}(S(u^{(1)}))(z)+1}{\bigl(\gamma_{1}
h_{q^{-1}}(S(u^{(1)}))(z)+q^{-1}z-\beta_{0}\bigr)\bigl(\gamma_{1}
S(u^{(1)})(z)+z-\beta_{0}\bigr)}
\\
\qquad{}
=
\frac{B(z)}{\bigl(\gamma_{1}
h_{q^{-1}}(S(u^{(1)}))(z)+q^{-1}z-\beta_{0}\bigr)\bigl(\gamma_{1}
S(u^{(1)})(z)+z-\beta_{0}\bigr)}\\
\qquad\quad{} -\frac{C(z)}{\bigl(\gamma_{1}
S(u^{(1)})(z)+z-\beta_{0}\bigr)}+D(z).
\end{gather*}
Equivalently
\begin{gather*}
\gamma_{1}
\bigl\{(h_{q^{-1}}\Phi)(z)+\big(q^{-1}-1\big)z\bigl(C(z)-(z-\beta_{0})D(z)\bigr)\bigr\}H_{q^{-1}}(S(u^{(1)}))(z)\\
\qquad{}=
\gamma_{1}^{2} D(z) S(u^{(1)})(z) h_{q
^{-1}}(S(u^{(1)}))(z)+\gamma_{1} \bigl\{((q^{-1}+1)z-2\beta_{0})
D(z)-C(z)\bigr\}S(u^{(1)})(z)
\\
\qquad\quad{} +B(z)+\big(q^{-1}z-\beta_{0}\big)(z-\beta_{0})D(z)
-\big(q^{-1}z-\beta_{0}\big)C(z)-(h_{q^{-1} }\Phi)(z).
\end{gather*}
Therefore the $q$-Riccati equation satisf\/ied by $S(u^{(1)})$
\begin{gather}
\big(h_{q ^{-1}}\Phi^{(1)} \big)H_{q ^{-1}}S(u^{(1)})=B^{(1)} S\big(u^{(1)}\big)h_{q
^{-1}}S\big(u^{(1)}\big)+C^{(1)} S\big(u^{(1)}\big) + D^{(1)} , \label{4.10}
\end{gather}
where
\begin{gather}
 K \Phi^{(1)}(x)=\Phi(x)+(q-1)x\{(qx-\beta_0)(h_qD)(x)-(h_qC)(x)\}  , \nonumber \\
 K B^{(1)}(x)=\gamma_1 D(x)  , \qquad
 K C^{(1)}(x )=\gamma_{1} \bigl\{\big(\big(q^{-1}+1\big)x-2\beta_0\big)D(x)-C(x)\bigr\}  ,\nonumber\\
 K D^{(1)}(x)=B(x)+(q^{-1}x-\beta_0)(x-\beta_0)D(x)
-(q^{-1}x-\beta_0)C(x)-(h_{q^{-1} }\Phi)(x) .
\label{4.11}
\end{gather}
$u^{(1)}$ is then a $q$-Laguerre--Hahn form.

Moreover, the regular form $u^{(1)}$ fulf\/ils the $q$-dif\/ference
equation
\begin{gather}
H_{q}\big(\Phi^{(1)} u^{(1)}\big)+\Psi^{(1)} u^{(1)}+B^{(1)}\big(x^{-1}u^{(1)}
h_q u^{(1)}\big)=0,\label{4.12}
\end{gather}
with
\begin{gather}
\Psi^{(1)}=-q^{-1}\big(C^{(1)}+ H_{q^{-1}}\Phi^{(1)}\big).\label{4.13}
\end{gather}
Likewise, it is straightforward to prove that the class of $u^{(1)}$
is also~$s$.
\end{proof}

\begin{example}\label{example2} If $u$ is a $q$-classical form satisfying the
$q$-analog of the distributional equation of Pearson type~\eqref{4.7} then
the associated~$u^{(1)}$ of~$u$ is a $q$-Laguerre--Hahn form of class
zero. $u^{(1)}$~and the formal Stieltjes function $S(u^{(1)})$
satisfy, respectively, the $q$-dif\/ference equation~\eqref{4.12} and the
$q$-Riccati equation~\eqref{4.10} where on account of~\eqref{4.11} and~\eqref{4.13}
\begin{gather*}
 K \Phi^{(1)}(x)=q   \frac{\phi''(0)}{2}   x^{2}+\left\{q
\phi'(0)+(q-1) \left(q \psi(0)+\beta_{0}\left(\frac{\phi''(0)}{2}+q
\psi'(0)\right)\right)\right\} x+\phi(0)  ,
\\
 K
\Psi^{(1)}(x)=-q^{-1}\left\{(q+1)\frac{\phi''(0)}{2}-\psi'(0))x+(q+1)\phi'(0) \right.\\
\left.\phantom{K\Psi^{(1)}(x)=}{} +q^{2}\psi(0)+
\big(q^{2}-q+2\big)\left(\frac{\phi''(0)}{2}+q\psi'(0)\right)\beta_{0}\right\} ,
\\
 K B^{(1)}(x)=-\gamma_{1}   \left(\frac{\phi''(0)}{2}+q
\psi'(0)\right) , \\
 K C^{(1)}(x )=\gamma_{1}   \left\{ -\psi'(0) x+\beta_{0}(\phi''(0)+2q \psi'(0))+q \psi(0)+\phi'(0)\right\}  ,\\
 K   D^{(1)}(x)=\psi(\beta_{0})
x-\phi(\beta_{0})-q\beta_{0}\psi(\beta_{0}).
\end{gather*}
\end{example}

\subsection[The inverse of a $q$-Laguerre-Hahn form]{The inverse of a $\boldsymbol{q}$-Laguerre--Hahn form}

Let $u$ be a regular form and $\{P_n\}_{n\geq0}$ its corresponding
MOPS satisfying~\eqref{1.1}. Let $\big\{P_{n}^{(1)}\big\}_{n\geq0}$ be its
associated sequence of the f\/irst kind fulf\/illing~\eqref{4.7} and
orthogonal with respect to the regular form $u^{(1)}$. The inverse
form of $u$ satisf\/ies~\cite[equation~(5.27)]{20}
\begin{gather}
x^{2} u^{-1}=-\gamma_{1} u^{(1)}. \label{4.14}
\end{gather}
The following results can be found in~\cite{2}
\begin{gather}
u^{-1}=\delta-\big(u^{-1}\big)_{1} \delta'-\gamma_{1} x^{-2} u^{(1)}.
\label{4.15}
\end{gather}
In general, the form $u^{-1}$ given by~\eqref{4.15} is regular if and only
if $\Delta_{n}\neq0$, $n\geq0$, with
\begin{gather*}
\Delta_{n}=\langle u^{(1)},\big(P_{n}^{(1)}\big)^{2}\rangle
\left\{\gamma_{1}+\sum_{\nu=0}^{n}
\frac{\bigl(\gamma_{1}P_{\nu-1}^{(2)}(0)-(u^{-1})_{1}P_{\nu}^{(1)}(0)\bigr)^{2}}{\langle
u^{(1)},(P_{\nu}^{(1)})^{2}\rangle}\right\} , \qquad n \geq 0,
\end{gather*}
where $\big\{P_{n}^{(2)}\big\}_{n\geq0}$ is the associated sequence of
$\big\{P_{n}^{(1)}\big\}_{n\geq0}$. In this case, the orthogonal sequence
$\big\{P_{n}^{(-)}\big\}_{n\geq0}$ relative to $u^{-1}$ is given by
\begin{gather*}
   P_{0}^{(-)}(x)=1 ,\qquad  P_{1}^{(-)}(x)=P_{1}^{(1)}(x)+b_{0} , \\
   P_{n+2}^{(-)}(x)=P_{n+2}^{(1)}(x)+b_{n+1} P_{n+1}^{(1)}(x)+a_{n} P_{n}^{(1)}(x),\qquad  n \geq 0 ,
  \end{gather*}
where
\begin{gather*}
b_{0}=\beta_{1}-\big(u^{-1}\big)_{1} , \\
b_{n+1}=\beta_{n+2}- \frac{\bigl((u^{-1})_{1}P_{n}^{(1)}(0)-\gamma_{1}
P_{n-1}^{(2)}(0)\bigr)
\bigl((u^{-1})_{1}P_{n+1}^{(1)}(0)-\gamma_{1}
P_{n}^{(2)}(0)\bigr)}{\Delta_{n}} , \qquad  n \geq 0 ,\\
a_{n}= \frac{\Delta_{n+1}}{\Delta_{n}} , \qquad  n \geq
0.
\end{gather*}
Also, the sequence $\big\{P_{n}^{(-)}\big\}_{n\geq0}$ satisf\/ies the
three-term recurrence relation
\begin{gather*}
   P_{0}^{(-)}(x)=1 ,\qquad  P_{1}^{(-)}(x)=x-\beta_{0}^{(-)} , \\
   P_{n+2}^{(-)}(x)=\big(x-\beta_{n+1}^{(-)}\big) P_{n+1}^{(-)}(x)-
   \gamma_{n+1}^{(-)} P_{n}^{(-)}(x),\qquad  n \geq 0,
  \end{gather*}
with
\begin{gather*}
\beta_{0}^{(-)}=\big(u^{-1}\big)_{1} ,\qquad
\beta_{n+1}^{(-)}=\beta_{n+2}+b_{n}-b_{n+1} ,
\qquad  n \geq 0,\\
\gamma_{1}^{(-)}=-\Delta_{0} , \qquad \gamma_{2}^{(-)}=\gamma_{1}
\frac{\Delta_{1}}{\Delta_{0}^{2}} , \qquad \gamma_{n+3}^{(-)}=
\frac{\Delta_{n+2} \Delta_{n}}{\Delta_{n+1}^{2}}
\gamma_{n+2},\qquad n \geq 0.
  \end{gather*}
In particular, when $\gamma_{1} >0$ and $u^{(1)}$ is positive
def\/inite, then $u^{-1}$ is regular.
When $u^{(1)}$ is symmetrical, then $u^{-1}$ is a symmetrical
regular form and we have
\begin{gather}
a_{2n}=\frac{\gamma_{1} \Lambda_{n}+1}{\gamma_{1} \Lambda_{n-1}+1}
  \gamma_{2n+2}  , \qquad a_{2n+1}=\gamma_{2n+3}  , \qquad
n\geq 0, \label{4.16}
\\
\gamma_{1}^{(-)}=-\gamma_{1} , \qquad \gamma_{2n+2}^{(-)}=a_{2n} ,
\qquad \gamma_{2n+3}^{(-)}=\frac{\gamma_{2n+2} \gamma_{2n+3}}{a_{2n}}
, \qquad n\geq 0, \label{4.17}
\end{gather}
with
\begin{gather}
\Lambda_{-1}=0 , \qquad \Lambda_{n}=\sum_{\nu=0}^{n}
\left(\prod_{k=0}^{\nu}\frac{\gamma_{2k+1}}{\gamma_{2k+2}}\right) ,
\quad n\geq 0 , \qquad \gamma_{0}=1.\label{4.18}
\end{gather}

\begin{proposition}\label{proposition11} If $u$ is a $q$-Laguerre--Hahn form
of class $s$, then, when $u^{-1}$ is regular,  $u^{-1}$ is a
$q$-Laguerre--Hahn form of class at most $s+2$.
\end{proposition}

\begin{proof} Let $u$ be a $q$-Laguerre--Hahn form
of class $s$ satisfying~\eqref{1.17}. It is seen in Proposition~\ref{proposition10} that
$u^{(1)}$ is also a $q$-Laguerre--Hahn form of class $s$ satisfying
the $q$-dif\/ference equation \eqref{4.12} with polynomials $\Phi^{(1)}$,
$\Psi^{(1)}$, $B^{(1)}$ respecting \eqref{4.11} and \eqref{4.13}.

Let us suppose $u^{-1}$ is regular that is to say $\Delta_{n}\neq0$,
$n\geq0$. Multiplying \eqref{4.12} by $(-\gamma_{1})$ and on account of
\eqref{4.14} and \eqref{1.7}, the $q$-dif\/ference equation \eqref{4.12} becomes
\begin{gather*}
H_{q}\big(x^{2}\Phi^{(1)}(x) u^{-1}\big)+x^{2}\Psi^{(1)}(x)
u^{-1}-q^{-2}\gamma_{1}^{-1}B^{(1)}\big(x^{-1}\big(x^{2}u^{-1}\big)
\big(x^{2}h_q u^{-1}\big)\big)=0.
\end{gather*}
Consequently, the form $u^{-1}$ satisf\/ies the following
$q$-dif\/ference equation
\begin{gather}
H_{q}\big(\Phi^{(-)} u^{-1}\big)+\Psi^{(-)} u^{-1}+B^{(-)}\big(x^{-1}u^{-1} h_q
u^{-1}\big)=0,\label{4.19}
\end{gather}
with
\begin{gather}
K \Phi^{(-)}(x)=x^{2} \big\{\Phi^{(1)}(x)+(1-q)\gamma_{1}^{-1} x (qx-\beta_0)(h_q B^{(1)})(x)\big\} ,  \nonumber\\
K\Psi^{(-)}(x)=x \Big\{\big(q^{-1}+1\big) \big(\big(h_{q^{-1}}
\Phi^{(1)}\big)(x)-q^{-1}
\Phi^{(1)}(x)\big)-q^{-3} x(H_{q^{-1}} \Phi^{(1)})(x)\nonumber\\
\phantom{K\Psi^{(-)}(x)=}{} +\gamma_{1}^{-1}
x\big(\big(2q^{-1}+q^{-2}-q^{-3}\big)x-\big(1+2q^{-2}-q^{-3}\big)\beta_{0}\big)B^{(1)}(x)\nonumber\\
\phantom{K\Psi^{(-)}(x)=}{}
 -\big(q^{-2}-1\big)\gamma_{1}^{-1}x(qx-\beta_{0})\big(h_q B^{(1)}\big)(x)\nonumber\\
\phantom{K\Psi^{(-)}(x)=}{}
-q^{-4} x^{2}(1-q)\gamma_{1}^{-1}(qx-\beta_{0})\big(H_q B^{(1)}\big)(x)-x C^{(1)}(x)\Big\} , \label{4.20}\\
 K B^{(-)}(x)=-\gamma_{1}^{-1} q^{-2} x^{4} B^{(1)}(x) .\tag*{\qed}
\end{gather}
\renewcommand{\qed}{}
\end{proof}

\begin{example} \label{example3}  Let $\mathcal{Y}(b,q^{2})$ be the form of Brenke type
which is symmetrical $q$-semiclassical of class one such that \cite[equation~(3.22), $q\leftarrow q^{2}$]{14}
\begin{gather}
H_{q}\big(x \mathcal{Y}\big(b,q^{2}\big)\big)-(b(q-1))^{-1}
\big(q^{-2}x^{2}+b-1\big) \mathcal{Y}\big(b,q^{2}\big)=0 \label{4.21}
\end{gather}
for $q\in \widetilde{\mathbb{C}}$, $ b\neq 0$, $b\neq q$, $
b\neq q^{-2n}$, $n\geq0$ and its MOPS $\{P_{n}\}_{n\geq0}$
satisfying \eqref{1.1} with~\cite{7}
\begin{gather}
    \beta_{n}=0 , \nonumber\\
    \gamma_{2n+1}=q^{2n+2}\bigl(1-b q^{2n}\bigr) , \qquad
    \gamma_{2n+2}=bq^{2n+2}\bigl(1-q^{2n+2}\bigr) , \qquad n\geq
    0 .\label{4.22}
\end{gather}
Denoting $\mathcal{Y}^{(1)}(b,q^{2})$ its associated form and
$\mathcal{Y}^{-1}(b,q^{2})$ its inverse one.
Taking into account~\eqref{4.21} we have
\begin{gather}
\Phi(x)=x, \qquad \Psi(x)=-(b(q-1))^{-1}
\big(q^{-2}x^{2}+b-1\big), \qquad B(x)=0.\label{4.23}
\end{gather}
Also, by virtue of \eqref{3.2} and \eqref{4.23} we get
\begin{gather}
C(x)=(b(q-1))^{-1} q^{-1} x^{2}+q (q-1)^{-1}\big(1-b^{-1}\big)-1,
\qquad D(x)=(bq(q-1))^{-1}x.\label{4.24}
\end{gather}
According to Proposition~\ref{proposition10} the form $\mathcal{Y}^{(1)}(b,q^{2})$
is $q$-Laguerre--Hahn of class one satisfying the $q$-dif\/ference
equation~\eqref{4.12} and its formal Stieltjes function satisf\/ies the
$q$-Riccati equation~\eqref{4.10} where on account of~\eqref{4.22}--\eqref{4.24} we
obtain for~\eqref{4.11}, \eqref{4.13}
\begin{gather}
    K \Phi^{(1)}(x)=b^{-1}x , \nonumber\\
    K \Psi^{(1)}(x)= -q^{-2}(b(q-1))^{-1} x^{2}+q(q-1)^{-1}\big(1-b^{-1}\big)-(qb)^{-1}-1, \nonumber\\
    K B^{(1)}(x)= \big(b^{-1}-1\big)q(q-1)^{-1} x, \nonumber\\
    K C^{(1)}(x)= q^{-2}(b(q-1))^{-1} x^{2}+1-q(q-1)^{-1}\big(1-b^{-1}\big), \nonumber\\
    K D^{(1)}(x)= q^{-2}(b(q-1))^{-1} x. \label{4.25}
\end{gather}
On the one hand, $\mathcal{Y}^{(1)}(b,q^{2})$ is a symmetrical
regular form, then $\mathcal{Y}^{-1}(b,q^{2})$ is also a symmetrical
regular form and we have for \eqref{4.16}--\eqref{4.18} according to \eqref{4.22}
\begin{gather*}
\Lambda_{-1}=0 ,\qquad \Lambda_{0}=\frac{b^{-1}-1}{1-q^{2}} , \qquad
\Lambda_{n}=\sum_{\nu=1}^{n+1} b^{-\nu}
\frac{(b;q^{2})_{\nu}}{(q^{2};q^{2})_{\nu}} , \qquad n\geq 1,
\\
\gamma_{1}^{(-)}=q^{2} (b-1), \qquad
\gamma_{2n+2}^{(-)}=bq^{2n+2}\bigl(1-q^{2n+2}\bigr)   \frac{1+
q^{2} (1-b)\Lambda_{n}}{1+ q^{2} (1-b) \Lambda_{n-1}}    , \qquad n
\geq 0,
\\
\gamma_{2n+3}^{(-)}=q^{2n+4} (1-b q^{2n+2})   \frac{1+ q^{2} (1-b)
\Lambda_{n-1}}{1+ q^{2} (1-b) \Lambda_{n}}   , \qquad n\geq 0,
\end{gather*}
with \cite{7}
\begin{gather*}
(a;q)_{0}=1, \qquad (a;q)_{n}=\prod_{k=1}^{n} \big(1-aq^{k-1}\big) ,
\qquad n \geq 1.
\end{gather*}
On the other hand, according to Proposition~\ref{proposition11}, \eqref{4.20} and \eqref{4.25},
the inverse form $\mathcal{Y}^{-1}(b,q^{2})$ is symmetrical
$q$-Laguerre--Hahn satisfying the $q$-dif\/ference equation~\eqref{4.19}
where
\begin{gather*}
K \Phi^{(-)}(x)=b^{-1}x^{3}\big(1-qx^{2}\big) ,\nonumber\\
K \Psi^{(-)}(x)=b^{-1}(q-1)^{-1} x^{2}
\big(b-q-q^{-3}(q-1)+\big({-}2q^{-4}+2q^{-3}+q^{-2}-q^{-1}+q\big)x^{2}\big),\nonumber\\
K B^{(-)}(x)=-b^{-1} q^{-3} (q-1)^{-1} x^{5}.
\end{gather*}
Thus, according to~\eqref{2.17} it is possible to simplify by $x$ one time
uniquely. Consequently, by virtue of~\eqref{2.16} the inverse form
$\mathcal{Y}^{-1}(b,q^{2})$ is $q$-Laguerre--Hahn of class two
fulf\/illing the $q$-dif\/ference equation
\begin{gather*}
H_{q}\big(x^{2}\big(x^{2}-q^{-1}\big)\mathcal{Y}^{-1}\big(b,q^{2}\big)\big)-q^{-1}x\big\{1+q(q-1)^{-1}\big(b-q-q^{-3}(q-1)\big)
\\
\qquad{}+\big(q(q-1)^{-1}\big({-}2q^{-4}+2q^{-3}+q^{-2}-q^{-1}+q\big)-q\big)
x^{2}\big\}\mathcal{Y}^{-1}\big(b,q^{2}\big) \\
\qquad{}+ q^{-3}(q-1)^{-1} x^{4}
\big(x^{-1}\mathcal{Y}^{-1}\big(b,q^{2}\big)h_{q}\mathcal{Y}^{-1}\big(b,q^{2}\big)\big)=0.
\end{gather*}
\end{example}

\subsection*{Acknowledgments}

 The authors are very grateful to the
referees for the constructive and valuable comments and
recommendations and for making us pay attention to a certain
references.

\pdfbookmark[1]{References}{ref}
\LastPageEnding

\end{document}